\input amstex
\documentstyle {amsppt}
\UseAMSsymbols \vsize 18cm \widestnumber\key{ZZZZZ}

\catcode`\@=11
\def\displaylinesno #1{\displ@y\halign{
\hbox to\displaywidth{$\@lign\hfil\displaystyle##\hfil$}&
\llap{$##$}\crcr#1\crcr}}
\def\ldisplaylinesno #1{\displ@y\halign{
\hbox to\displaywidth{$\@lign\hfil\displaystyle##\hfil$}&
\kern-\displaywidth\rlap{$##$} \tabskip\displaywidth\crcr#1\crcr}}
\catcode`\@=12

\refstyle{A}

\let \ol=\overline
\let \ul=\underline
\let \ti=\widetilde

\font\nr=eufb7 at 10pt

 \font\srm=cmr10 at 7.5pt

\font\main=cmsy10 at 10pt

\font\smain=cmsy10 at 7.5pt \font\ssmain=cmsy10 at 5.625pt

\font \fin=lasy8 at 15.4 pt
\def \X{\mathop{\hbox{\nr X}^{\hbox{\srm nr}}}\nolimits}

\def \o{\mathop{\hbox{\main O}}\nolimits}
\def \so{\mathop{\hbox{\smain O}}\nolimits}
\def \sso{\mathop{\hbox{\ssmain O}}\nolimits}

\def \H{\mathop{\hbox{\main H}}\nolimits}

\def \End{\mathop{\hbox{\rm End}}\nolimits}

\def \Hom{\mathop{\hbox{\rm Hom}}\nolimits}

\def \GL{\mathop{\hbox{\rm GL}}\nolimits}

\topmatter
\title Alg\`ebres de Hecke avec
param\`etres et Repr\'esentations d'un groupe $p$-adique
classique: Pr\'eservation du spectre temp\'er\'e\endtitle

\rightheadtext{Pr\'eservation du spectre temp\'er\'e}
\author Volker Heiermann \endauthor
\address Aix-Marseille Univ,
IML, 13288 Marseille C\'edex 9, France;
CNRS, UMR 6206, 13288 Marseille C\'edex 9, France \endaddress

\email volker.heiermann\@univmed.fr \endemail

\thanks Ce travail a b\'enifici\'e d'une aide de l'Agence
Nationale de la Recherche portant la r\'eference
ANR-08-BLAN-0259-02.
\endthanks

\abstract  Let $G$ be an orthogonal or symplectic $p$-adic group
(not necessarily split) or an inner form of a general linear
$p$-adic group. In a previous paper, it was shown that the
Bernstein components of the category of smooth representations of
$G$ are equivalent to the category of right modules over some
Hecke algebra with parameters, or more general over the
semi-direct product of such an algebra with a finite group
algebra.

The aim of the present paper is to show that this equivalence preserves the tempered spectrum and the discrete series representations.

\null\noindent{R\'ESUM\'E: } Soit $G$ un groupe $p$-adique orthogonal ou symplectique (pas n\'ecessaire-ment d\'eploy\'e) ou une forme int\'erieure d'un groupe lin\'eaire g\'en\'eral $p$-adique. Dans un article pr\'ec\'edent, il a \'et\'e montr\'e que les composantes de Bernstein de la cat\'egorie des repr\'esentations lisses de $G$ sont \'equivalentes \`a la cat\'egorie des modules \`a droite sur une certaine alg\`ebre de Hecke avec param\`etres, ou plus g\'en\'eralement sur le produit semi-direct d'une telle alg\`ebre avec l'alg\`ebre d'un groupe fini.

Le but de l'article pr\'esent est de montrer que cette \'equivalence pr\'eserve le spectre temp\'er\'e ainsi que la s\'erie discr\`ete.
\endabstract

\endtopmatter
\document
Fixons un corps local non archim\'edien $F$ de valeur absolue
normalis\'ee $\vert\cdot\vert _F$. Soit $G$ (le groupe des points
rationnels d') un groupe orthogonal ou symplectique ou encore
d'une forme int\'erieure d'un groupe lin\'eaire g\'en\'eral
(d\'efini sur $F$). On consid\'erera la cat\'egorie $Rep(G)$ des
repr\'esentations complexes lisses de $G$.

Soient $P$ un sous-groupe parabolique de $G$, $M$ un sous-groupe
de Levi de $P$ et $(\sigma ,E)$ une repr\'esentation
irr\'eductible cuspidale de $M$. Notons $\o _{\sigma }$ ou encore
$\o $ l'ensemble des tordus de $\sigma $ par un caract\`ere non
ramifi\'e de $M$, modulo \'equivalence, appel\'e \it classe
inertielle \rm de $\sigma $, et $Rep_{\so }(G)$ la composante de
Bernstein de $Rep(G)$ associ\'ee \`a $\o $. \'Ecrivons $M^1$ pour
l'intersection des noyaux des caract\`eres non ramifi\'es de $M$,
$ind_{M^1}^M$ pour l'induction compacte et $i_P^G$ pour le
foncteur de l'induction parabolique normalis\'ee. Fixons une
composante irr\'eductible $(\sigma _1,E_1)$ de la restriction de
$\sigma $ \`a $M^1$. Posons $E_{B_{\sso }}=ind_{M^1}^ME_1$ et
$\H=End_G(i_P^GE_{B_{\sso }})$.

Il a \'et\'e remarqu\'e par J. Bernstein \cite{Ro} que
$i_P^GE_{B_{\sso }}$ est un prog\'en\'erateur de $Rep_{\so }(G)$.
Par le contexte de Morita qui en r\'esulte, on a une \'equivalence
entre $Rep_{\so }(G)$ et la cat\'egorie des $\H-\hbox{modules \`a
droite}$, donn\'ee par $V\mapsto \Hom_G(i_P^GE_{B_{\sso }},$ $V)$.

Dans \cite{H1}, l'auteur de cet article a montr\'e que $\H $ est
le produit semi-direct d'une alg\`ebre de Hecke avec param\`etres
(au sens de Lusztig \cite{L1}) par l'alg\`ebre d'un groupe fini.
L'objet de l'article pr\'esent est de prouver que l'\'equivalence
de cat\'egories $V\mapsto \Hom_G(i_P^GE_{B_{\sso }},V)$ pr\'eserve
le spectre temp\'er\'e ainsi que les repr\'esentations de carr\'e
int\'egrable (modulo le centre de $G$).

L'article est organis\'e de la mani\`ere suivante: au premier
paragraphe, nous rassemblons les r\'esultats n\'ecessaires
relatifs aux repr\'esentations temp\'er\'ees et de carr\'e
int\'egrable des groupes $p$-adiques. Le deuxi\`eme paragraphe est
consacr\'e \`a celles des alg\`ebres de Hecke avec param\`etres.
Au troisi\`eme paragraphe, nous prouvons un lien entre les
exposants pour les modules d'alg\`ebres de Hecke avec param\`etres
et ceux des modules de Jacquet de repr\'esentations de $G$.
L'objet du quatri\`eme paragraphe est d'\'enoncer certaines
propri\'et\'es des modules de Jacquets des repr\'esentations de
$G$ dont nous aurons besoin au paragraphe 5 pour prouver notre
th\'eor\`eme principal.

Je remercie A.-M. Aubert pour la correction de quelques erreurs de frappe et E. Opdam pour des discussions sur la d\'efinition des repr\'esentations de carr\'e int\'egrable (modulo le centre) des alg\`ebres de Hecke. 

\null{\bf 1.} Fixons un tore d\'eploy\'e maximal $A_0$ dans $M$.
Notons $W$ le groupe de Weyl de $G$ associ\'e \`a ce tore. Un
sous-groupe parabolique $P'$ de $G$ sera dit semi-standard s'il
contient $A_0$. Il existe alors un unique sous-groupe de Levi $M'$
de $P'$ qui contient $A_0$. On dira que $M'$ est un sous-groupe de
Levi semi-standard (de $P'$ ou de $G$).  On \'ecrira $A_{M'}$ pour
le tore maximal d\'eploy\'e contenu dans le centre de $M'$,
$a_{M'}$ pour son alg\`ebre de Lie r\'eelle, $a_{M'}^*$ pour son
dual et $a_{M',\Bbb C}^*$ pour le complexifi\'e de  $a_{M'}^*$. On
a une d\'ecomposition canonique $a_{M'}^*=a_{M'}^{G*}\oplus
a_G^*$. Le groupe de Weyl de $M'$ relatif \`a $A_0$ sera not\'e
$W^{M'}$.

Il existe un unique sous-groupe parabolique semi-standard minimal
$P_0$ contenu dans $P$. On dira que $P'$ et $M'$ sont standard
s'ils contiennent respectivement $P_0$ et le sous-groupe de Levi
semi-standard, not\'e $M_0$, de $P_0$.

Si $P$ est un sous-groupe parabolique semi-standard de sous-groupe
de Levi standard $M$, on notera $\ol{P}$ le sous-groupe
parabolique semi-standard de $G$ qui v\'erifie $\ol{P}\cap P=M$.
On dira qu'un sous-groupe parabolique semi-standard $P$ de $G$ est
\it anti-standard \rm si $\ol{P}$ est standard.

Lorsque $M'$ est un sous-groupe de Levi semi-standard, on notera $\X (M')$ le groupe de ses caract\`eres non ramifi\'es. On a une surjection canonique $a_{M',\Bbb C}^*\rightarrow \X(M')$ que l'on notera $\lambda\mapsto\chi_{\lambda }$. Sa restriction \`a  $a_{M'}^*$ est injective d'image \'egale au sous-groupe des caract\`eres \`a valeurs r\'elles. On d\'efinit $\Re(\chi _{\lambda })=\Re(\lambda )$. Remarquons que, lorsque $\chi $ est un caract\`ere de $M$, le caract\`ere $\vert\chi\vert _F$ est non ramifi\'e \`a valeurs r\'eelles. On peut donc d\'efinir la partie r\'eelle de $\chi $, not\'ee $\Re(\chi )$, comme \'etant celle de  $\vert\chi\vert _F$.

Soit $(\pi ,V)$ une repr\'esentation admissible de $M'$.  Appelons sous-espace caract\'eristique la donn\'ee d'un caract\`ere non ramifi\'e $\chi $ de $A_{M'}$ et d'un sous-espace non nul $V_{\chi }$ form\'e des \'el\'ements $v$ de $V$ pour lesquels il existe un entier $N\geq 0$ tel que, pour tout $a'\in A_{M'}$, $$(\pi (a')-\chi (a'))^Nv=0.$$ On appelle $\chi $ un exposant du $A_{M'}$-module $V$, et on sait que $V$ est, en tant que $A_{M'}$-module, la somme directe de ces sous-espaces caract\'eristiques.

Le foncteur adjoint \`a gauche de $i_{P'}^G$ sera not\'e
$r_{P'}^G$. Rappelons que, par le deuxi\`eme th\'eor\`eme
d'adjonction de J. Bernstein, $r_{\ol{P'}}^G$ est l'adjoint \`a
droite de $i_{P'}^G$.

\null {\bf 1.1} Lorsque $P'$ est un sous-groupe parabolique
standard de $G$ de sous-groupe de Levi standard $M'$, notons
$\Sigma (P')$ l'ensemble des restrictions non nulles \`a $A_{M'}$
des racines r\'eduites de $A_0$ dans l'alg\`ebre de Lie du radical
unipotent de $P'$, $\Delta (P')$ le sous-ensemble des restrictions
non nulles des racines simples (relatives \`a $P_0$), et
$\ol{^+a_{P'}^{G*}}$ (resp. $^+a_{P'}^{G*}$) l'ensemble des
combinaisons lin\'eaires de $\Delta (P')$ \`a coefficients $\geq
0$ (resp. $>0$).

Pour les repr\'esentations admissibles de $G$, on a le crit\`ere suivant, d\^u \`a Casselman:

\null{\bf Proposition:} \it (\cite{W, III.1.1 et III.2.2}) Soit $\pi$ une repr\'esentation admissible de $G$.

(i) Pour que la repr\'esentation $\pi $ soit temp\'er\'ee, il faut et il suffit que, pour tout sous-groupe parabolique standard $P'$ de $G$ tel que $P'=G$ ou $P'$ soit propre et maximal, et tout exposant $\chi $ de $r_{P'}^G\pi $, on ait $\Re (\chi )\in\ \ol{^+a_{P'}^{G*}}$.

(ii) Supposons que $\pi $ admette un caract\`ere central. Pour que $\pi $ soit de carr\'e int\'egrable (modulo le centre de $G$), il faut et il suffit que le caract\`ere central soit unitaire et que, pour tout sous-groupe parabolique standard $P'$ de $G$, propre et maximal, et tout exposant $\chi $ de $r_{P'}^G\pi $, on ait $\Re (\chi )\in\ ^+a_{P'}^{G*}$. \rm

\null {\bf 1.2} Remarquons que le groupe de Weyl $W$ agit sur l'ensemble des sous-groupe de Levi semi-standard de $G$.
On d\'esignera cette action par $M'\mapsto wM'$. Lorsque $M'$ est un sous-groupe de Levi
semi-standard, notons $^{M'}W^M$ l'ensemble des doubles-classes
$W^{M'}\backslash W/W^M$ et $^{M'}W^M(M')$ le sous-ensemble des
doubles-classes des \'el\'e-ments $w$ avec $wM\subseteq M'$.
Le groupe $^{M'}W^{M'}(M')$ (ainsi que les repr\'esentants de ses
\'el\'ements dans $W$) agit sur l'ensemble des sous-groupes
paraboliques semi-standard de $G$ de sous-groupe de Levi $M'$. On
notera cette action $P'\mapsto wP'$.

\null{\bf Proposition:} \it Soit $\pi $ une repr\'esentation irr\'eductible dans $Rep_{\so }(G)$. Alors,

(i) la repr\'esentation $\pi $ est temp\'er\'ee, si et seulement si les parties r\'eelles des exposants de $r_P^G\pi $ sont dans $\ol{^+a_P^{G*}}$.

(ii) la repr\'esentation $\pi $ est  de carr\'e int\'egrable
(modulo le centre de $G$), si et seulement si son caract\`ere
central est unitaire et si les parties r\'eelles des exposants de
$r_P^G\pi $ sont dans $^+a_P^{G*}$.

\null Preuve: \rm Il r\'esulte de la transitivit\'e du foncteur de Jacquet que les exposants de $r_P^G\pi $ ont cette propri\'et\'e, si cette propri\'et\'e est v\'erifi\'ee par rapport \`a tout sous-groupe parabolique maximal (cf. \cite{W, III.1.1 et III.2.2}). La r\'eciproque tient au fait que $\pi $ se plonge dans l'induite parabolique d'une repr\'esentation cuspidale d'un sous-groupe de Levi de $P$. Explicitons la preuve:

Fixons un \'el\'ement $\sigma '$ de $\o $ tel que $\pi $ soit une
sous-repr\'esentation de $i_P^G\sigma '$. Soit $P'$ un sous-groupe
parabolique standard maximal de $G$ de sous-groupe de Levi
standard $M'$. Suite \`a la cuspidalit\'e de $\sigma '$, $r_{M\cap
w^{-1}P'}^{M}\sigma '$ est nulle si $M\not\subseteq w^{-1}M'w$ et
\'egal \`a $\sigma '$ sinon. La semi-simplifi\'ee de
$r_{P'}^Gi_P^G\sigma '$ est donc, par le lemme g\'eom\'etrique
\cite{BZ2, 2.12}, \'egale \`a $$\bigoplus _{w\in\ ^{M'}W^M(M')}
i_{M'\cap wP}^{M'}w\sigma '.$$ Pour des raisons analogues, la
semi-simplifi\'ee de $r_P^Gi_P^G\sigma '$ est \'egale \`a
$$\bigoplus _{w\in\ ^{M}W^M(M)}w\sigma'.$$

Notons $W_{P'}(\pi )$ (resp. $W_P(\pi )$) l'ensemble des $w\in\ ^{M'}W^M(M')$ (resp. $w\in\ ^{M}W^M(M)$) tels que $i_{M'\cap wP}^{M'}w\sigma'$ et $r_{P'}^G(\pi )$ ont un sous-quotient en commun (resp. $w\sigma '$ est un sous-quotient de $r_P^G\pi $).

Supposons les parties r\'eelles des exposants de $r_P^G\pi $ dans $\ol{^+a_P^{G*}}$ (resp. $^+a_P^{G*}$). On a donc $\Re (\chi _{w\sigma '})=\Re (w\chi _{\sigma '})\in\
\ol{^+a_P^{G*}}$ (resp. $^+a_P^{G*}$) pour $w\in W_P(\pi )$.

La semi-simplifi\'ee de $r_{P\cap M'}^{M'}i_{M'\cap
wP}^{M'}w\sigma '$ est $$\bigoplus _{w'\in\
^{M}(W^{M'})^{wM}(M)}w'w\sigma '.$$ Si $w\in W_{P'}(\pi )$, alors,
par la transitivit\'e du foncteur de Jacquet, $w'w\in W_P(\pi )$
pour au moins un $w'\in\ ^{M}(W^{M'})^{wM}(M)$, et alors $\Re
(w'w\chi _{\sigma '})\in \ol{^+a_P^{G*}}$ (resp. $^+a_P^{G*}$). Le
caract\`ere central de $i_{M'\cap wP}^{M'}w\sigma'$ est $(\chi
_{w\sigma '})_{\vert A_{M'}}=(w'w\chi _{\sigma '})_{\vert
A_{M'}}$. C'est dans $\ol{^+a_{P'}^{G*}}$ (resp. $^+a_{P'}^{G*}$).

Signalons finalement que, lorsque $(\pi ,V)$ est temp\'er\'ee, le fait que le caract\`ere central $\chi _{\pi }$ de $\pi $ v\'erifie $\Re(\chi _{\pi })=0$ r\'esulte du fait que $\chi _{\pi }$ est la restriction au centre de $G$ du caract\`ere central $\chi _{\sigma '}$ de $\sigma '$ et que $\Re(\chi _{\sigma '})\in \ol{^+a_P^{G*}}$, remarquant que $\sigma '$ est un sous-quotient de $r_P^G\pi $ par r\'eciprocit\'e de Frob\'enius. \hfill{\fin 2}

\null{\bf 2.} Consid\'erons maintenant le cas d'une alg\`ebre
complexe $\H $ qui est le produit semi-direct d'une alg\`ebre de
Hecke avec param\`etres par l'alg\`ebre d'un groupe fini. Une
alg\`ebre de Hecke avec param\`etres (cf. \cite{L}) est associ\'ee
\`a un quintuplet $\Xi =(\Lambda ,\Sigma ,\Lambda^{\vee },\Sigma
^{\vee },\Delta )$ avec $\Lambda $, $\Lambda ^{\vee }$ des
r\'eseaux en dualit\'e, $\Sigma $ et $\Sigma ^{\vee }$ des
syst\`emes de racines dans $\Lambda $ et $\Lambda ^{\vee }$
respectivement, duaux l'un de l'autre, et $\Delta $ une base de
$\Sigma $, ainsi que \`a un syst\`eme de param\`etres $q_{\alpha
},q_{\alpha }'$, $\alpha\in\Delta $.  Nous noterons $\H(\Xi
,\{q_{\alpha }, q_{\alpha }'\})$ l'alg\`ebre de Hecke associ\'ee
et $R$ le groupe fini tel que $\H =\H(\Xi ,\{q_{\alpha },
q_{\alpha }'\})\rtimes \Bbb C[R]$. On supposera que $R$ agisse sur
$\Lambda $. (Remarquons que cette hypoth\`ese est bien
v\'erifi\'ee dans \cite{H1}, mais fausse pour un groupe r\'eductif
g\'en\'eral (cf. \cite{GR} pour un exemple).)

Soit $V$ un $\H $-module de dimension finie. Appelons sous-espace caract\'eristique de $V$ la donn\'ee d'un caract\`ere non ramifi\'e $\chi $ de $\Lambda $ et d'un sous-espace non nul $V_{\chi }$ form\'e des \'el\'ements $v$ de $V$ pour lesquels il existe un entier $N\geq 0$ tel que, pour tout $\lambda \in\Lambda $, $$(\lambda -\chi (\lambda ))^Nv=0.$$ On appellera $\chi $ un exposant de $V$. On sait que $V$ est, en tant que $\Lambda $-module, la somme directe de ces sous-espaces caract\'eristiques.

Nous dirons que $\vert\chi\vert $ est la valeur absolue de $\chi $ et faisons alors la d\'efinition suivante:

\null{\bf 2.1 D\'efinition:} \it Soit $V$ un $\H $-module \`a droite de dimension finie.

(i) On dit que $V$ est temp\'er\'e, si les valeurs absolues des
exposants de $V$ pour l'action de $\Lambda $ sont des combinaisons
lin\'eaires de $\Delta ^{\vee }$ \`a coefficients $\leq 0$.

(ii) Soit $Z$ un sous-r\'eseau de $\Lambda $, contenu dans le
centre de $\H $. On dit que $V$ est de carr\'e int\'egrable modulo
$Z$, si le caract\`ere pour l'action de $Z$ est unitaire,  si
$Z\cup \Delta$ engendre un $\Bbb Z$-module d'indice fini de
$\Lambda $ et si les valeurs absolues des exposants de $V$ pour
l'action de $\Lambda $ sont des combinaisons lin\'eaires de
$\Delta ^{\vee }$ \`a coefficients $<0$.\rm

\null\it Remarque: \rm Il suit de (ii) qu'une condition
n\'ecessaire pour qu'un $\H $-module $V$ poss\`ede des
repr\'esentations de carr\'e int\'egrable modulo $Z$ est que
$Z\cup \Delta $ engendre un sous-module d'indice fini du $\Bbb
Z$-module $\Lambda $.

Cette d\'efinition est bien compatible avec celles donn\'ees dans
\cite{KL}, \cite{O}, \cite{DO} et \cite{BC}

\null{\bf 3.} Notons $M^{\sigma }$ l'ensemble des \'el\'ements $m$
de $M$ tels que $^m\sigma _1\simeq\sigma _1$. Rappelons que c'est
un sous-groupe d'indice fini de $M$ \cite{H1, 1.16}. L'alg\`ebre
de Hecke affine contenue dans $\H=\End _G(i_P^GE_{B_{\sso }})$ est
associ\'ee au r\'eseau $\Lambda _{\so }:=M^{\sigma }/M^1$ et \`a un
syst\`eme de racines $\Sigma _{\so }$ que l'on pr\'ecisera plus
tard (cf. section {\bf 4.}). Soulignons ici seulement que les
\'el\'ements de $\Sigma _{\so }^{\vee }$ sont des multiples \`a
coefficients $>0$ de certaines racines dans $\Sigma (P)$ (via le
plongement de $\Lambda _{\so }^{\vee }$ dans $a_M^*$).

Nous noterons $b_m$ l'image d'un \'el\'ement $m$ de $M^{\sigma }$ dans $\Lambda _{\so }$ (ou plus g\'en\'eralement dans $M/M^1$) et $Z_G$ l'image du tore central d\'eploy\'e $A_G$ dans $\Lambda _{\so }$.

Dans les notations de \cite{H1}, on a $\Bbb C[\Lambda _{\so }]=B_{\so }$ et donc $B_{\so }^{\times }=\Lambda _{\so }$. Par cons\'equent, on \'ecrira dans la suite $B_{\so }^{\times }$ \`a la place de $\Lambda _{\so }$.

Posons $B=\Bbb C[M/M^1]$. La repr\'esentation $(ind_{M^1}^M\sigma _{\vert M^1}, ind_{M^1}^ME _{\vert M^1})$ est isomorphe \`a la repr\'esentation $\sigma _B$ de $M$ dans l'espace $E_B:=E\otimes _{\Bbb C}B$, d\'efinie par $\sigma _B(m)e\otimes b=\sigma(m)e\otimes bb_m$. La repr\'esentation $(\sigma _{B_{\sso }}, E_{B_{\sso }})$ est isomorphe \`a une sous-repr\'esentation de $(\sigma _B, E_B)$. L'espace $E_{B_{\sso }}$ est muni d'une structure de $B_{\so }$-module, induite par la structure de $B$-module de $E_B$, apr\`es restriction des scalaires.

\null{\bf 3.1 } Notons $\chi _{\sigma }$ le caract\`ere central de $\sigma $.

\null {\bf Lemme:} \it Soit $a\in A_G$. Alors, pour tout $v\in i_P^GE_{B_{\sso }}$, on a $$b_av=\chi _{\sigma }(a)^{-1}(i_P^G\sigma _{B_{\sso }})(a)v.$$

\null Preuve: \rm Il suffit de montrer que ceci est vrai pour $\sigma _B$ et $v\in i_P^GE_B$. Soit $g\in G$. Alors, $v(g)=\sum _ie_i\otimes b_i$ avec $e_i\in E_i$, $b_i\in B$, et, donc
$$\eqalign{((i_P^G\sigma _B)(a)v)(g)&=v(ga)=v(ag)=\sigma _B(a)v(g)\cr &=\sum _i\sigma (a)e_i\otimes b_ib_a=\chi _{\sigma }(a)\sum _ie_i\otimes b_ib_a\cr &=\chi _{\sigma }(a)b_av(g).\cr }$$ Ceci prouve le lemme. \hfill{\fin 2}

\null{\bf 3.2 Proposition:} \it Soit $(\pi ,V)$ une repr\'esentation irr\'eductible de $G$ qui appartient \`a $Rep_{\so }(G)$. Notons $\chi _{\pi }$ le caract\`ere central de $\pi $. Alors, $Z_G$ agit par un caract\`ere sur $\Hom _G(i_P^GE_{B_{\sso }},V)$ donn\'e par  $z\mapsto (\chi _{\pi }\chi _{\sigma }^{-1})(z)$.

En particulier, l'action de $Z_G$ se fait par un caract\`ere unitaire, si et seulement si $\chi _{\pi }$ est un caract\`ere unitaire.

\null Preuve: \rm Soit $\varphi $ dans $\Hom _G(i_P^GE_{B_{\sso
}},V)$, $v$ dans $i_P^GE_{B_{\sso }}$ et $a$ dans $A_G$. Alors,
par le lemme {\bf 3.1}, $$\varphi (b_av)=\varphi (\chi _{\sigma
}(a)^{-1}(i_P^G\sigma _B)(a)v)=\chi _{\sigma }(a)^{-1}\pi
(a)\varphi (v)=\chi _{\sigma }(a)^{-1}\chi _{\pi }(a)\varphi
(v).$$

Choisissons un caract\`ere non ramifi\'e $\chi _{\lambda }$ de $M$
tel que $\pi $ se plonge dans la repr\'esentation induite
$i_P^G(\sigma\otimes\chi_{\lambda })$. Alors, $\chi _{\pi }$ est
la restriction de $\chi _{\sigma }\chi _{\lambda }$ au centre de
$G$. La restriction de $\chi _{\sigma }^{-1}\chi _{\pi }$ \`a
$A_G$ est donc bien trivial sur $A_G\cap M^1$, i.e. l'action de
$A_G$ se factorise bien par $Z_G$. \hfill{\fin 2}

\null {\bf 3.3} Soit $M'$ un sous-groupe de Levi contenant $M$.  Remarquons que $A_{M'}\subseteq A_M\subseteq M^{\sigma }$. Notons $B_{A_{M'}}^{\times }$ le sous-groupe de $B_{\so }^{\times }$ form\'e des \'el\'ements $\{b_{a'}\vert a'\in A_{M'}\}$. Tout $B_{\so }^{\times }$-module $V$ de dimension finie se d\'ecompose en sous-espaces caract\'eristiques pour l'action de $B_{A_{M'}}^{\times }$.

\null{\bf Proposition:} \it Soit $(\pi ,V)$ une repr\'esentation irr\'eductible lisse de $G$ qui appartient \`a $Rep _{\so }(G)$. Soit $P'$ un sous-groupe parabolique standard de $G$ qui contient $P$, de sous-groupe de Levi standard $M'$.

Alors, $\Hom _G(i_P^GE_{B_{\sso }},V)$ se d\'ecompose en
sous-espaces caract\'eristiques pour l'ac-tion de
$B_{A_{M'}}^{\times }$. Notons $Exp _{P'}(\pi )$ l'ensemble des
exposants de $r_{\ol{P'}}^G\pi $. Les exposants de $\Hom
_G(i_P^GE_{B_{\sso }},V)$ sont les caract\`eres $\chi\chi_{\sigma
}^{-1}$ avec $\chi\in Exp_{P'}(\pi )$, et v\'erifiant $\Hom _{M'}$
$(i_{P\cap M'}^{M'}E_{B_{\sso }},(r_{\ol{P'}}^GV )_{\chi })$ $\ne
0$.

\null Preuve: \rm Par le deuxi\`eme th\'eor\`eme d'adjonction de Bernstein (cf. \cite{B, III.3} et \cite{R, VI.9}), on a un isomorphisme d'espaces vectoriels
$$\Hom _G(i_P^GE_{B_{\sso }},V)=\Hom _{M'}(i_{P\cap M'}^{M'}E_{B_{\sso }},r_{\ol{P'}}^GV).\eqno{\hbox{\rm (*)}}$$
D\'ecomposons d'abord $\Hom _{M'}(i_{P\cap M'}^{M'}E_{B_{\sso
}},r_{\ol{P'}}^GV)$ en espaces caract\'eristiques pour l'action de
$B_{A_{M'}}$. Pour $\chi\in Exp_{P'}(\pi )$, notons
$\Hom_{M'}(i_{P\cap M'}^{M'}E_{B_{\sso }},r_{\ol{P'}}^GV)_{\chi }$
le sous-espace de $\Hom _{M'}($ $i_{P\cap M'}^{M'}E_{B_{\sso }},$
$r_{\ol{P'}}^GV)$ form\'e des homomorphismes dont l'image est
contenue dans $(r_{\ol{P'}}^G\pi )_{\chi }$. On a $$\Hom
_{M'}(i_{P\cap M'}^{M'}E_{B_{\sso }},r_{\ol{P'}}^GV)=\bigoplus
_{\chi\in Exp_{P'}(\pi )}\Hom _{M'}(i_{P\cap M'}^{M'}E_{B_{\sso
}},r_{\ol{P'}}^GV)_{\chi }.$$ Fixons $\chi\in Exp_{P'}(\pi )$.
Comme $r_{\ol{P'}}^G\pi $ est de longueur finie, il existe un
entier $N_{\chi }$ tel que, pour tout $a'\in A_{M'}$ et $v\in
(r_{\ol{P'}}^GV)_{\chi }$,
$$((r_{\ol{P'}}^G\pi )(a')-\chi (a'))^{N_{\chi }}v=0.$$ Soit
$\varphi\in \Hom _{M'}(i_{P\cap M'}^{M'}E_{B_{\sso
}},r_{\ol{P'}}^GV)_{\chi }$, et prouvons que $\varphi (b_{a'}-\chi
(a'))^{N_{\chi }}=0$ pour tout $a'\in A_{M'}$.

Par le lemme {\bf 3.1}, on a, pour $a'\in A_{M'}$, $v'\in i_{P\cap
M'}^{M'}E_{B_{\sso }}$,
$$\eqalign {\varphi (b_{a'}v')&=\varphi(\chi _{\sigma
}(a')^{-1}(i_{P\cap M'}^{M'}\sigma _{B_{\sso }})(a')v')\cr
&=\chi _{\sigma }^{-1}(a')(r_{\ol{P'}}^G\pi)(a')\varphi(v').\cr }$$

Il en suit par la formule du bin\^ome que
$$\eqalign{&\varphi((b_{a'}-(\chi_{\sigma }^{-1}\chi) (a'))^{N_{\chi }}v')\cr
=&\varphi(\sum _{k=0}^{N_{\chi }}{\binom {N_{\chi }}
k}b_{a'}^k(-\chi _{\sigma }^{-1}(a')\chi (a'))^{{N_{\chi
}}-k}v')\cr =&\sum _{k=0}^{N_{\chi }}\binom {N_{\chi }} k (\chi
_{\sigma }^{-1}(a')(r_{\ol{P'}}^G\pi )(a'))^k(-\chi _{\sigma
}^{-1}(a')\chi (a'))^{{N_{\chi }}-k}\varphi(v')\cr =&\chi _{\sigma
}^{-{N_{\chi }}}(a')((r_{\ol{P'}}^G\pi )(a')-\chi (a'))^{N_{\chi
}}\varphi (v')\cr =&0.\cr }$$ Les exposants dans $\Hom
_{M'}(i_{P\cap M'}^{M'}E_{B_{\sso }},r_{\ol{P'}}^GV)$ pour
l'action de $B_{A_{M'}}^{\times }$ sont donc les caract\`eres
$\chi\chi _{\sigma }^{-1}$, $\chi\in Exp _{P'}(\pi )$, tels que
$\Hom _{M'}(E_{B_{\sso }},(r_{\ol{P'}}^GV)_{\chi })\ne 0$.

Traduisons ces exposants en exposants de $\Hom _G(i_P^GE_{B_{\sso
}},V)$ via l'isomorphisme (*). Le deuxi\`eme th\'eor\`eme
d'adjonction de Bernstein est r\'ealis\'e par un  isomorphisme
$$\iota :\Hom _G(i_{P'}^Gi_{P \cap M'}^{M'}E_{B_{\sso
}},V)\rightarrow \Hom _{M'}(i_{P\cap M'}^{M'}E_{B_{\sso
}},r_{\ol{P'}}^GV),$$ d\'efini de la mani\`ere suivante (cf. \cite{B, III.3.1}: notons
$r_{\ol{P'}}^{G,P'\ol{P'}}$ la restriction de $r_{\ol{P'}}^G$ au
sous-espace $(i_{P'}^Gi_{P \cap M'}^{M'}E_{B_{\sso
}})_{P'\ol{P'}}$ form\'e des \'el\'ements $v$ avec
$supp(v)\subseteq P'\ol{P'}$. Cette restriction est injective, et
l'image de $r_{\ol{P'}}^{G,P'\ol{P'}}$ est une
sous-repr\'esentation de $r_{\ol{P'}}^Gi_{P'}^G i_{P \cap
M'}^{M'}E_{B_{\sso }}$ qui est isomorphe \`a $i_{P\cap
M'}^{M'}E_{B_{\sso }}$. Cet isomorphisme est compatible avec
l'action de $B_{\so }$. Identifions ces deux espaces au moyen d'un
tel isomorphisme.

Avec cette identification, l'application $\iota $ associe \`a un
homomorphisme $\varphi: i_{P'}^G$ $i_{P \cap M'}^{M'}E_{B_{\sso
}}\rightarrow V$ l'homomorphisme $\iota(\varphi ):i_{P\cap
M'}^{M'}$ $E_{B_{\sso }}\rightarrow r_{\ol{P'}}^GV$ qui envoie un
\'el\'ement $r_{\ol{P'}}^Gv$ de l'image de
$r_{\ol{P'}}^{G,P'\ol{P'}}$ sur $(r_{\ol{P'}}^G\varphi
)(r_{\ol{P'}}^Gv):=r_{\ol{P'}}^G(\varphi (v)).$ Montrons que
$\iota $ est compatible avec l'action de $B_{A_{M'}}^{\times }$.

Remarquons d'abord que, pour $r_{\ol{P'}}^Gv$ dans l'image de
$r_{\ol{P'}}^{G,P'\ol{P'}}$ et $a'\in A_{M'}$, on a, par le lemme
{\bf 3.1} et la $B_{A_{M'}}^{\times }$-lin\'earit\'e de
$r_{\ol{P'}}^{G}$,
$$\eqalign{r_{\ol{P'}}^{G}v &=b_{a'}^{-1}\chi_{\sigma
}(a')^{-1}(i_{P\cap M'}^{M'}\sigma _{B_{\sso }})(a')r_{\ol{P'}}^{G}v\cr
&=r_{\ol{P'}}^{G}(b_{a'}^{-1}\chi _{\sigma
}(a')^{-1}(i_{P'}^Gi_{P\cap M'}^{M'}\sigma _{B_{\sso
}})(a')v).\cr}$$

Par suite, $$\eqalign{&\iota(\varphi b_{a'})(r_{\ol{P'}}^{G}v) \cr
=&r_{\ol{P'}}^G(\varphi (\chi_{\sigma }(a')^{-1}i_{P'}^Gi_{P\cap
M'}^{M'}\sigma_{B_{\sso }}(a')v))\cr =&\iota(\varphi
)r_{\ol{P'}}^{G}(\chi _{\sigma }(a')^{-1}i_{P'}^Gi_{P\cap
M'}^{M'}\sigma _{B_{\sso }}(a')v)\cr =&\iota(\varphi )\chi
_{\sigma }(a')^{-1}i_{P\cap M'}^{M'}\sigma _{B_{\sso
}}(a')r_{\ol{P'}}^{G}(v)\cr =&\iota(\varphi
)b_{a'}r_{\ol{P'}}^{G}(v),\cr }$$ d'o\`u la $A_{M'}$-lin\'earit\'e
de $\iota $.

Remarquons que $\iota^{-1}(\Hom_{M'}(i_{P\cap M'}^{M'}E_{B_{\sso
}}, r_{\ol{P'}}^GV)_{\chi })$ est l'ensemble des $\varphi $ tels
que $r_{\ol{P'}}^G(\varphi((i_{P'}^Gi_{P\cap M'}^{M'}E_{B_{\sso
}})_{P'\ol{P'}}))\subseteq (r_{\ol{P'}}^GV)_{\chi }$. Pour
$\varphi $ dans cet espace et $a'\in A_{M'}$, on a donc avec la
formule du bin\^ome
$$0=\iota(\varphi )(b_{a'}-(\chi_{\sigma }\chi )(a'))^{N_{\chi
}}=\iota (\varphi (b_{a'}-(\chi_{\sigma }\chi )(a'))^{N_{\chi
}}),$$ d'o\`u $(b_{a'}-(\chi_{\sigma }\chi )(a'))^{N_{\chi
}}\varphi =0,$ puisque $\iota $ est injective et donc $$\iota(\Hom
_G(i_{P'}^Gi_{P \cap M'}^{M'}E_{B_{\sso }},V)_{\chi\chi_{\sigma
}^{-1}})=\Hom _{M'}(i_{P\cap M'}^{M'}E_{B_{\sso
}},r_{\ol{P'}}^GV)_{\chi\chi_{\sigma }^{-1}}.$$

Pour obtenir l'isomorphisme (*), il faut composer avec
l'isomorphisme $$\kappa :i_{P'}^Gi_{P\cap M'}^{M'}E_{B_{\sso
}}\rightarrow i_P^GE_{B_{\sso }}$$ qui envoie un \'el\'ement $v$
de l'espace \`a gauche sur l'application $G\rightarrow E_{B_{\sso
}}$ donn\'ee par $g\mapsto (v(g))(1)$. Plus pr\'ecis\'ement,
l'isomorphisme (*) envoie un \'el\'ement $\psi $ de $\Hom
_G(i_P^GE_{B_{\sso }},V)$ sur $\iota (\psi\circ\kappa )$. Comme
$\kappa $ est $B_{\so }$-lin\'eaire, il en est ainsi pour
l'application $\psi\mapsto\psi\circ\kappa $. L'isomorphisme (*)
laisse donc invariant les espaces caract\'eristiques, et
$\Hom_G(i_P^GE_{B_{\sso }}, V)_{\chi\chi _{\sigma }^{-1}}$ est
l'ensemble des homomorphismes $\varphi $ v\'erifiant
$r_{\ol{P'}}^G(\varphi((i_P^GE_{B_{\sso }})_{P'\ol{P'}}))\subseteq
(r_{\ol{P'}}^GV)_{\chi }$, ce qui \'equivaut bien \`a $\Hom _{M'}$
$(i_{P\cap M'}^{M'}$ $E_{B_{\sso }},(r_{\ol{P'}}^GV )_{\chi })$
$\ne 0$.\hfill{\fin 2}

\null{\bf 3.4 Corollaire:} \it L'ensemble des parties r\'eelles
des exposants du $B_{\so }^{\times }$-module $\Hom
_G(i_P^GE_{B_{\sso }},$ $ V)$ est \'egal \`a celui  des
caract\`eres centraux des quotients irr\'eductibles de
$r_{\ol{P}}^GV$ qui appartiennent \`a $\o $.

\null Preuve: \rm Remarquons d'abord que tout sous-quotient
irr\'eductible de $r_{\ol{P}}^GV$ est isomorphe \`a un quotient
irr\'eductible, les composantes de $r_{\ol{P}}^GV$ \'etant des
repr\'esenta-tions cuspidales de $M$. Ainsi, les exposants de
$r_{\ol{P}}^GV$ s'identifie aux caract\`eres centraux des
quotients irr\'eductibles de $r_{\ol{P}}^GV$.

Le $B_{\so }^{\times }$-module $\Hom _G(i_P^GE_{B_{\sso }}, V)$ se
d\'ecompose en sous-espaces caract\'eristiques. Les restrictions
\`a $B_{A_M}^{\times }$ des exposants de ces sous-espaces
caract\'eris-tiques sont des exposants pour $\Hom
_G(i_P^GE_{B_{\so }}, V)$ en tant que $B_{A_M}^{\times }$-module.
Donc, tout exposant pour l'action de $B_{A_M}^{\times }$ est la
restriction d'un exposant pour l'action de $B_{\so }^{\times }$.
Rappelons que les caract\`eres de $B_{\so }^{\times }$ sont des
caract\`eres non ramifi\'es de $M$ et que ceux de $B_{A_M}^{\times
}$ sont des caract\`eres non ramifi\'es pour $A_M$. Or, la
restriction de $M$ \`a $A_M$ est injective sur les caract\`eres
non ramifi\'es \`a valeurs r\'eelles. Les parties r\'eelles des
exposants pour l'action par $B_{\so }^{\times }$ ou par
$B_{A_M}^{\times }$ sont donc les m\^emes. Par la proposition {\bf
3.3}, celles-ci correspondent aux parties r\'eelles des
caract\`eres non ramifi\'es $\chi $ tels que $\Hom _M (E_{B_{\sso
}},(r_{\ol{P}}^G\pi )_{\chi\chi _{\sigma }})\ne 0$. Comme $\chi
_{\sigma }$ est un caract\`ere unitaire, le corollaire en
r\'esulte. \hfill{\fin 2}

\null{\bf 4.} Il r\'esulte presque directement de la proposition
{\bf 1.2} et du corollaire {\bf 3.4} que les spectres discrets et
temp\'er\'es sont pr\'eserv\'es par l'\'equivalence de cat\'egorie
si les \'el\'ements de $\Delta_{\so }^{\vee }$ sont les multiples
des \'el\'ements de $\Delta (P)$ par une constante non nulle.
Cette hypoth\`ese n'est toutefois pas toujours v\'erifi\'ee. Pour
\'etablir le cas g\'en\'eral, il faut du travail suppl\'ementaire.
On va d'abord collecter quelques propri\'et\'es qui nous serons
utile pour le calcul du module de Jacquet.

\null {\bf 4.1} \it D\'efinition: \rm On dira que $\o $ est \it en position standard \rm si

\null (i) le sous-groupe de Levi $\underline{M}$ de groupe des points $F$-rationnels \'egal \`a $M$ est \'egal \`a

$$\GL_{k_1}\times\cdots\times\GL _{k_1}
\times\GL_{k_2}\times\cdots \times\GL _{k_2}\times
\cdots\times\GL_{k_r}\times\cdots\times\GL_{k_r}\times\underline{H}_k,$$

o\`u $\underline{H}_k$ d\'esignant un groupe
semi-simple de rang absolu $k$ et du m\^eme type que $G$;

\null et

\null (ii) $\o $ contient une repr\'esentation unitaire de la
forme

$$\sigma =\sigma_1\otimes\cdots\otimes\sigma
_1\otimes\sigma_2\otimes\cdots\otimes \sigma
_2\otimes\cdots\otimes\sigma_r\otimes\cdots\otimes\sigma
_r\otimes\tau, $$

les classes inertielles des $\sigma _i$ \'etant deux \`a deux
distinctes avec $\sigma _i\simeq\sigma _i^{\vee }$ si $\sigma _i$
et $\sigma _i^{\vee }$ sont dans une m\^eme orbite inertielle,
$\sigma _i^{\vee }\not\in\o_{\sigma _j}$ lorsque $i\ne j$.

\null On appelera $\sigma $ un \it point de base \rm de $\o $, si,
en outre, pour tout $i$, ou bien il existe $s>0$ tel que l'induite
parabolique de $\sigma _i\vert\cdot\vert _F^s$ soit r\'eductible
ou bien l'induite parabolique de $\sigma _i\vert\cdot\vert _F^s$
est irr\'eductible pour tout nombre complexe $s$ de partie
r\'eelle $\ne 0$.

\null \it Remarque: \rm On a vu dans \cite{H1} (en particulier
\cite{H1, 1.13}), que l'on peut toujours se ramener \`a $\o $
standard et alors trouver un point de base dans $\o $, v\'erifiant
les hypoth\`eses de la d\'efinition {\bf 4.1} sauf celle que
$\sigma _i^{\vee }$ n'est pas dans la classe inertielle de $\sigma
_j$ pour $i\ne j$. Or, il est facile d'arriver \`a cette condition
suppl\'ementaire apr\`es conjugaison par un \'el\'ement du groupe
de Weyl.

\null {\bf 4.2} Supposons $\o $ en position standard et d\'esignons, pour $i=1,\dots ,r$, par $d_i$ le nombre des facteurs \'egal \`a $\sigma _i$ dans la d\'ecomposition de $\sigma $. Notons $a=(a_{1,1}, \dots , a_{1,d_1}, a_{2,1}\dots , a_{r,d_r},$ $1)$ les \'el\'ements du centre d\'eploy\'e $A_M$ de $M$ et $\alpha _{i,j}$ les \'el\'ements de $\Delta (P)$ avec  $$\alpha _{i,j}(a)=\cases a_{i,j}a_{i,j+1}^{-1} & \hbox{\rm si } j\ne d_i, \cr a_{i,d_i}a_{i+1,1}^{-1} & \hbox{\rm si } j=d_i\ \hbox{\rm et } i\ne r, \cr a_{r,d_r} & \hbox{\rm sinon}. \cr\endcases$$

\null{\bf Proposition:} (\cite{H1, 1.13 et  6.1} ) \it  Fixons un point de base $\sigma $ de $\o $. Les composantes irr\'eductibles $\Sigma _{\so ,i}$ de $\Sigma _{\so }$ correspondent aux facteurs $\sigma _i$ de $\sigma $.

La composante irr\'eductible $\Sigma _{\so ,i}$ est

(i) de type $B_{d_i}$ s'il existe $s>0$ telle que l'induite parabolique de $\sigma _i\vert\cdot\vert^s\otimes\tau $ soit r\'eductible;

(ii) de type $D_{d_i}$ si les conclusions de (i) ne sont pas v\'erifi\'ees, mais $\sigma _i\simeq\sigma _i^{\vee }$.

(iii) de type $A_{d_i-1}$ dans les autre cas.

Une base $\Delta _{\so ,i}^{\vee }$ de l'ensemble $\Sigma _{\so
,i}^{\vee }$ des racines duales est alors form\'ee respectivement
des \'el\'ements des ensembles $\{\alpha _{i,1},\dots ,\alpha
_{i,d_i-1},\alpha _{i,d_i}+\cdots +\alpha _{r,d_r}\}$, $\{\alpha
_{i,1},\dots ,$ $\alpha _{i,d_i-1},\alpha _{i,d_i-1}+2\alpha
_{i,d_i}\}$ et $\{\alpha _{i,1},\dots ,\alpha _{i,d_i-1}\}$,
multipli\'e chacun par un nombre r\'eel $>0$ convenable.

Lorsque $\Sigma _{\so ,i}$ est une composante
irr\'eductible de $\Sigma _{\so }$ de type $D_{d_i}$, notons $r_i$
la sym\'etrie simple associ\'ee \`a la racine courte du syst\`eme
de racine de type $B_{d_i}$ qui contient $\Sigma _{\so ,i}^{\vee
}$. Alors, $r_i$ laisse invariant la base $\Delta _{\so
,i}^{\vee }$.

Le groupe $R$ associ\'e \`a $Rep _{\so }(G)$ est le groupe $R_{\so
}$ engendr\'e par les $r_i$ correspondant aux composantes
irr\'educibles de type $D_{d_i}$ de $\Sigma _{\so }$.

\null Remarque: \rm En fait, il est indispensable pour cette
proposition que le point de base $\sigma $ poss\`ede la condition
suppl\'ementaire de {\bf 4.1} qui manque dans \cite{H1}. Il y a
donc un erratum \`a \cite{H1}.

\null{\bf 4.3 Lemme:} \it Soit $\sigma'\in\o $ et soit $w\in W$
tel que $w\sigma '$ soit un sous-quotient irr\'eductible de
$r_P^G\pi $. Soit $\alpha $ un \'el\'ement de $\Delta (P)$ dont
aucun multiple n'appartient \`a $\Delta _{w\so }^{\vee }$ et
notons $\ti{P}_{\alpha }$ le sous-groupe parabolique standard de
$G$ de sous-groupe de Levi $s_{\alpha }M$.

Alors, $s_{\alpha }w\sigma '$ est un sous-quotient irr\'eductible
de $r_{\ti{P}_{\alpha }}^G\pi $.

\null Preuve: \rm Comme $w\sigma '$ est un sous-quotient
irr\'eductible de $r_P^G\pi $, c'est m\^eme un quotient de
$r_P^G\pi $. Rappelons que, lorsque $\alpha\in\Delta (P)$, la
composante connexe du centralisateur du sous-groupe de $A_M$
\'egal au noyau de $\alpha $ est un sous-groupe de Levi standard
de $G$ que l'on notera $M_{\alpha }$.

Par r\'eciprocit\'e de Frobenius, on a
$$0\ne \Hom _G(\pi ,i_P^Gw\sigma ')=\Hom _{M_{\alpha }}(r_{P_{\alpha }}^G\pi ,i_{P\cap M_{\alpha }}^{M_{\alpha }}w\sigma ').$$
Par les r\'esultats de Harish-Chandra \cite{H1, 1.2 et 1.8},
$i_{P\cap M_{\alpha }}^{M_{\alpha }}w\sigma '$ est ou bien une
repr\'esentation irr\'eductible ou bien $s_{\alpha }w\sigma'\simeq
w\sigma'$ et $s_{\alpha }P\cap M=\ol{P}\cap M_{\alpha }$.

Dans le premier cas, $s_{\alpha }w\sigma '$ est un sous-quotient
de $r_{\ti{P}_{\alpha }\cap M_{\alpha }}^{M_{\alpha }}i_{P\cap
M_{\alpha }}^{M_{\alpha }}w\sigma '$, et donc, par exactitude et
transitivit\'e du foncteur de Jacquet, de $r_{\ti{P}_{\alpha
}}^G\pi $.

Dans le deuxi\`eme cas, la conclusion de la proposition est
triviale, puisque ${\ti{P}_{\alpha }}=P$ et $s_{\alpha
}w\sigma'\simeq w\sigma '$. \hfill{\fin 2}

\null{\bf 4.4 Proposition:}  \it Supposons $\o $ en position
standard. Soit $(\pi ,V)$ une repr\'esen-tation irr\'eductible de
$G$ qui appartient \`a $Rep_{\so }(G)$ et soit $\sigma '$ un
sous-quotient irr\'eductible de $r_P^G\pi $.

i) Notons $w_{i,j}$ l'\'el\'element du groupe de Weyl qui permute
les facteurs $GL_{k_i}\times\cdots\times GL_{k_i}$ et
$GL_{k_j}\times\cdots\times GL_{k_j}$ de $M$ et $P_{i,j}$ le
sous-groupe parabolique standard de sous-groupe de Levi
$w_{i,j}Mw_{i,j}^{-1}$. Alors, $w_{i,j}\sigma '$ est un
sous-quotient irr\' eductible de $r_{P_{i,j}}^G\pi $.

ii) Si $\Sigma _{\so ,i}$ est de type $D_{d_i}$, alors $r_i\sigma
'$ est un sous-quotient irr\' eductible de $r_P^G\pi $ (o\`u $r_i$
est la sym\'etrie d\'efinie dans {\bf 4.2}).

iii) Supposons $\Sigma _{\so ,i}$ de type $A_{d_i-1}$ et notons
$w_{0,i}$ l'\'el\'ement le plus long de son groupe de Weyl
compos\'e avec l'\'el\'ement le plus long du groupe de Weyl du
syst\`eme de type $B_{d_i}$ qui contient $\Sigma _{\so ,i}$.
Alors, $w_{0,i}\sigma '$ est encore un sous-quotient
irr\'eductible de $r_P^G\pi $.

iv) Plus g\'en\'eralement, si $r$ est un \'el\'ement du groupe de
Weyl qui permute les facteurs de $\underline{M}$ et qui v\'erifie
$r\Delta_{\so }=\Delta_{r\so }$, alors $r\sigma '$ est  un
sous-quotient irr\' eductible de $r_{rP}^G\pi $.

\null Preuve: \rm On obtient (i) avec le lemme {\bf 4.3}, en
permutant successivement les facteurs adjoint de $\sigma '$ issus
d'orbites inertielles disjointes.

Pour (ii), on peut se ramener, gr\^ace \`a (i), \`a $i=r$. Alors,
$r_r$ provient d'une sym\'etrie simple et le lemme {\bf 4.3}
s'applique directement.

Pour (iii), observons que $r_i\sigma '$ change la composante
$\sigma _{i,d_i}'$ de $\sigma '$ en ${\sigma '_{i,d_i}}^{\vee }$
et que, par hypoth\`ese, $\sigma _{i,d_i}'$ et ${\sigma
_{i,d_i}'}^{\vee }$ ne se trouvent pas dans une m\^eme orbite
inertielle. On peut donc \'echanger ${\sigma _{i,d_i}'}^{\vee }$
avec les $\sigma _{i,j}$ qui le pr\'ec\`ede suite au lemme {\bf
4.3}. En r\'ep\'etant l'op\'eration, on aboutit \`a $w_{0,i}\sigma
'$.

Concernant (iv), la condition $r\Delta_{\so }=\Delta_{r\so }$
assure que $r$ est produit de sym\'etrie simple $s_{\alpha }$
v\'erifiant les hypoth\`ese du lemme {\bf 4.3}. \hfill{\fin 2}

\null{\bf 4.5} Notons $W(M)$ l'ensemble des \'el\'ements du groupe
de Weyl de $G$ qui stabilisent $M$ et $W(M,\o )$ le sous-ensemble
form\'e des \'el\'ements qui stabilisent en outre $\o $.

\null{\bf Lemme:} \it Soit $w\in W(M)$ tel que $w\sigma '$ soit un
sous-quotient irr\'eductible de $r_P^G\pi $. Alors, il existe
$w'\in W(M,\o )$ et $r\in W(M)$ tel que $w'\sigma '$ soit un
sous-quotient irr\'eductible de  $r_P^G\pi $, $r\Delta_{\so }=
\Delta_{r\so }$ et $w=rw'$.

\null Preuve: \rm En permutant successivement les facteurs de
$w\sigma '$ issus d'orbites inertielles disjointes par des
sym\'etries simples, on transforme $w\sigma '$ en un \'el\'ement
de $\o $. La proposition {\bf 4.4} (iv) s'applique \`a ces
sym\'etries. Notons $r$ le produit de ces sym\'etries successives.
Alors, par la proposition {\bf 4.4} (iv), $rw\sigma '$ est bien un
sous-quotient irr\'eductible de $r_P^G\pi $. Par ailleurs, par
construction, $r^{-1}\Delta _{\so }=\Delta _{r^{-1}\so
}$.

\hfill{\fin 2}

\null{\bf 5. Th\'eor\`eme:} \it Soit $(\pi ,V)$ une
repr\'esentation irr\'eductible de $G$ qui appartient \`a $Rep_{\so }(G)$.

Pour que la repr\'esentation $(\pi ,V)$ soit temp\'er\'ee (resp.
de carr\'e int\'egrable modulo le centre de $G$), il faut et il
suffit que le $\End _G(i_P^GE_{B_{\so}})$-module $\Hom
_G(i_P^GE_{B_{\so}}, V)$ soit temp\'er\'e (resp. de carr\'e
int\'egrable modulo $Z_G$).

\null\it Preuve: \rm Supposons le
$\End_G(i_P^GE_{B_{\sso}})$-module $\Hom _G(i_P^GE_{B_{\sso}}, V)$
de carr\'e int\'egra-ble modulo $Z_G$. En passant aux sous-groupes
paraboliques anti-standard, on d\'eduit de la proposition {\bf
1.2} qu'il faut montrer que les exposants des sous-quotients
irr\'educ-tibles de $r_{\ol{P}}^G\pi $ sont des combinaisons
lin\'eaires \`a coefficients $<0$ de $\Delta (P)$.

Par d\'efinition (cf. {\bf 2.1}), $\Delta _{\so }$ engendre
n\'ecessairement $M^{\sigma }\cap G^1/M^1$. Par suite, $a_M^{G*}$
est engendr\'e par $\Delta _{\so }^{\vee}$. Comme $\Delta _{\so
}^{\vee}$ est form\'e d'\'el\'ements qui sont des multiples \`a
coefficients $>0$ de $\Sigma(P)$, tout \'el\'ement de $\Delta
_{\so }^{\vee }$ est une combinaison lin\'eaire \`a coefficients
$\geq 0$ de $\Delta (P)$. Par notre hypoth\`ese, les parties
r\'eelles des exposants de $\Hom _G(i_P^GE_{B_{\so}}, V)$ relatifs
\`a l'action de $B_M$ sont des combinaisons lin\'eaires de $\Delta
_{\so }^{\vee }$ \`a coefficients $<0$. Il en est donc par le
corollaire {\bf 3.4} de m\^eme pour les parties r\'eelles des
exposants des sous-quotients irr\'eductibles de $r_{\ol{P}}^G\pi $
qui sont dans $\o$. Or, comme $\Delta _{\so }^{\vee }$ engendre
$a_M^{G*}$, celles-ci sont alors des combinaisons lin\'eaires \`a
coefficients $<0$ de $\Delta (P)$.

En g\'en\'eral, un sous-quotient irr\'eductible de
$r_{\ol{P}}^G\pi $ est de la forme $w\sigma'$ avec $\sigma '\in\o
$ et $w\in W(M)$. Par le lemme {\bf 4.3}, on peut \'ecrire $w=rw'$
avec $w'\in W(M,\o )$, $r\in W(M)$, tel que $w'\sigma '$ soit un
quotient irr\'eductible de $r_{\ol{P}}^G\pi $ et $r\Delta_{\so
}=\Delta _{r\so }$. Il en r\'esulte que l'exposant de $w\sigma '$
est combinaison lin\'eaire \`a coefficients $<0$ de $\Delta _{r\so
}^{\vee }$ et on en d\'eduit comme avant, en rempla\c cant $\Delta
_{\so }$ par $\Delta _{r\so }$, que le caract\`ere central de
$w\sigma '$ est une combinaison lin\'eaire \`a coefficient $<0$ de
$\Delta (P)$. Par cons\'equence, $\pi $ est bien de carr\'e
int\'egrable modulo le centre.

On prouve de mani\`ere analogue que $\Hom _G(i_P^GE_{B_{\so}}, V)$
temp\'er\'e implique $(\pi ,V)$ temp\'er\'ee. (Dans ce cas,
l'hypoth\`ese que $\Delta _{\so }$ engendre $M^{\sigma }\cap
G^1/M^1$ n'est pas n\'ecessaire.)

R\'eciproquement, supposons $(\pi ,V)$ de carr\'e int\'egrable
modulo le centre.  Alors, en passant aux sous-groupes paraboliques
anti-standard, on d\'eduit de la proposition {\bf 1.2} que les
caract\`eres centraux des sous-quotients irr\'eductibles de
$r_{\ol{P}}^G\pi $ sont des combinaisons lin\'eaires de $\Delta
(P)$ \`a coefficients $<0$.

Sous ces hypoth\`eses, $\Delta _{\so }^{\vee }$ engendre
$a_M^{G*}$ \cite{Si} (voir aussi \cite{H3}) et, pour des raisons
de dimension,  $\Delta _{\so }$ engendre $M^{\sigma }\cap G^1/M^1$.  Par
cons\'equence, $\Delta (P)$ et $\Delta _{\so }^{\vee }$ ont m\^eme
nombre d'\'el\'ements et engendrent tous les deux $a_M^{G*}$. On
en d\'eduit que $\Sigma _{\so }$ est ou bien irr\'eductible de
type $A_n$ (pour $\ul{G}$ une forme int\'erieure de $GL _n$) ou
bien ses composantes sont toutes de type $B_{d_i}$ ou $D_{d_i}$ .
Si $\ul{G}$ est une forme int\'erieure de $GL_n$, $\Delta_{\so
}^{\vee }$ est, modulo une constante $>0$, \'egal \`a $\Delta (P)$
et l'implication imm\'ediate par {\bf 3.4}.

Sinon, supposons d'abord que $\Sigma _{\so }$ soit produit de
syst\`emes de type $B_{d_i}$. Soit $\lambda $ la partie r\'eelle
du caract\`ere central d'un sous-quotient de $r_{\ol{P}}^G\pi $
qui appartient \`a $\o $. Par hypoth\`ese, $\lambda $ est
combinaison lin\' eaire de $\Delta (P)$ \`a coefficients $<0$, et
il faut montrer que $\lambda $ est combinaison lin\'eaire \`a
coefficients $<0$ de $\Delta _{\so }^{\vee }$. Rappelons (cf. {\bf
4.2}) que nous avons not\'e les \'el\'ements de $\Delta (P)$ par
$\alpha _{i,j}$, seule la racine $\alpha _{r,d_r}$ \'etant courte.
Notons les \'el\'ements de $\Delta _{\so }^{\vee }$ par $\alpha
_{i,j}'$. Apr\`es les avoir multipli\'es par un nombre r\'eel $<0$
convenable, on peut, suite \`a la proposition {\bf 4.2}, supposer
que  $\alpha _{i,j}'=\alpha _{i,j}$ pour $j\ne d_i$ et $\alpha
_{i,d_i}'=\alpha _{i,d_i}+\alpha _{i+1,1}+\cdots+\alpha _{r,d_r}$.
\'Ecrivant $\lambda =\lambda _{1,1}\alpha _{1,1}+\cdots +\lambda
_{r,d_r}\alpha _{r,d_r},$ on trouve $$\eqalign{\lambda &=\lambda
_{1,1}\alpha _{1,1}'+\cdots +\lambda _{1,d_1}\alpha _{1,d_1}'\cr
&+(\lambda _{2,1}-\lambda _{1,d_1})\alpha _{2,1}'+(\lambda
_{2,2}-\lambda _{1,d_1})\alpha _{2,2}'+\cdots +(\lambda
_{2,d_2}-\lambda _{1,d_1})\alpha _{2,d_2}'\cr &+(\lambda
_{3,1}-\lambda _{2,d_2})\alpha _{3,1}'+(\lambda _{3,2}-\lambda
_{2,d_2})\alpha _{3,2}'+\cdots +(\lambda _{3,d_3}-\lambda
_{2,d_2})\alpha _{3,d_3}'\cdots \cr &+(\lambda _{r,1}-\lambda
_{r-1,d_{r-1}})\alpha _{r,1}'+(\lambda _{r,2}-\lambda
_{r-1,d_{r-1}})\alpha _{r,2}'+\cdots +(\lambda _{r,d_r}-\lambda
_{r-1,d_{r-1}})\alpha _{r,d_r}'\cr}$$ Il faut donc montrer que
$\lambda _{i,j}-\lambda _{i-1,d_{i-1}}<0$ pour $i=2,\dots ,r$.
Notons $w_{1,i}$ l'\'el\'ement du groupe de Weyl qui permute les
facteurs $GL_{k_1}\times\cdots\times GL_{k_1}$ et
$GL_{k_i}\times\cdots\times GL_{k_i}$ de $\underline{M}$.
D'apr\`es le proposition {\bf 4.4} (ii), $w_{1,i}\lambda $ est
e\'galement la partie r\'elle du caract\`ere central d'un
sous-quotient irr\'eductible de $r_{\ol{P}}^G\pi $. Notons
$\ti{P}$ le sous-groupe parabolique standard de sous-groupe de
Levi $w_{1,i}M$ et $\ti{\alpha }_{i,j}$ les \'el\'ements de
$\Delta(\ti{P})$. On constate que
$$\eqalign{w_{1,i}\lambda=&\lambda _{1,1}\ti{\alpha }_{i,1}+\cdots +\lambda _{1,d_1-1}\ti{\alpha }_{i,d_1-1}+
\lambda _{1,d_1}(-\ti{\alpha }_{2,1}-\cdots -\ti{\alpha }_{i,d_i-1})\cr
&+\lambda _{2,1}\ti{\alpha }_{2,1}+\lambda _{2,2}\ti{\alpha }_{2,2}+\cdots +\lambda _{2,d_2}\ti{\alpha }_{2,d_2}+\cdots \cr
&+\lambda _{i-1,1}\ti{\alpha }_{i-1,1}+\cdots +\lambda _{i-1,d_{i-1}-1}\ti{\alpha }_{i-1,d_{i-1}-1}\cr
&+\lambda _{i-1,d_{i-1}}(-\ti{\alpha }_{1,1}-\cdots -\ti{\alpha }_{i-1,d_{i-1}-1})\cr
&+\lambda _{i,1}\ti{\alpha } _{1,1}+\cdots +\lambda _{i,d_i-1}\ti{\alpha }_{1,d_i-1}+\lambda _{i,d_i}(\ti{\alpha }_{1,d_i}+\cdots +\ti{\alpha }_{i,d_1})\cr
&+\lambda _{i+1,1}\ti{\alpha }_{i+1,1}+\cdots +\lambda _{r,d_r}\ti{\alpha }_{r,d_r}\cr
=&(\lambda _{i,1}-\lambda _{i-1,d_{i-1}})\ti{\alpha } _{1,1}+\cdots +(\lambda _{i,d_i-1}-\lambda _{i-1,d_{i-1}})\ti{\alpha }_{1,d_i-1}\cr &+(\lambda _{i,d_i}-\lambda _{i-1,d_{i-1}})\ti{\alpha }_{1,d_i}
+\cdots, \cr}$$ d'o\`u l'on d\'eduit que $\lambda _{i,j}-\lambda _{i-1,d_{i-1}}<0$ pour $j=1,\dots ,d_i$. Ceci prouve l'implication lorsque toutes les composantes irr\'eductibles de $\Sigma _{\so }$ sont de type $B_{d_i}$.

Supposons maintenant qu'il y a des composantes irr\'eductibles
$\Sigma _{\so ,i}$ de type $D_{d_i}$. On a $\Delta  _{\so
,i}=\{\alpha _{i,1}',\dots ,\alpha _{i,d_i-1}', \alpha
_{i,d_i-1}'+2\alpha _{i,d_i}'\},$
$$\eqalign{\lambda =&\lambda_{1,1}\alpha _{1,1}'+\cdots \lambda_{i-1,d_{i-1}}\alpha _{i-1,d_{i-1}}'\cr
&+\lambda_{i,1}\alpha _{i,1}' +\cdots +(\lambda
_{i,d_i-1}-{\lambda _{i,d_i}\over 2})\alpha _{i,d_{i-1}}'+{\lambda
_{i,d_i}\over 2}(\alpha _{i,d_{i-1}}'+2\alpha _{i,d_i}')\cr
&+\lambda_{i+1,1}\alpha _{i+1,1}'+\dots +\lambda_{r,d_r}\alpha
_{r,d_r}'\cr}$$ avec $\lambda _{i,j}<0$, et il faut donc montrer
que $\lambda _{i,d_i-1}-{\lambda _{i,d_i}\over 2}<0$.

En \'echangeant de nouveau les facteurs
$GL_{k_1}\times\cdots\times GL_{k_1}$ et
$GL_{k_i}\times\cdots\times GL_{k_i}$ de $\underline{M}$ (ce qui
ne revient maintenant qu'\`a permuter les $\alpha _{i_j}'$), on se
ram\`ene \`a $i=1$. D'apr\`es le proposition {\bf 4.4} (ii),
$r_1\lambda $ est encore combinaison lin\'eaire \`a coefficients
$<0$ de $\Delta (P)$. Or, on a
$$\eqalign{r_1\lambda
=&\lambda_{1,1}\alpha _{1,1}' +\cdots +\lambda _{1,d_1-1}(\alpha
_{1,d_1-1}'+2\alpha _{1,d_1}')-\lambda _{1,d_1}\alpha
_{1,d_1}'\cr &+\lambda_{2,1}\alpha _{2,1}'+\dots
+\lambda_{r,d_r}\alpha _{r,d_r}'\cr =& \lambda_{1,1}\alpha _{1,1}'
+\cdots +\lambda _{1,d_1-1}\alpha _{1,d_1-1}'+(2\lambda
_{1,d_1-1}-\lambda _{1,d_1})\alpha _{1,d_1}'\cr
&+\lambda_{2,1}\alpha _{2,1}'+\dots +\lambda_{r,d_r}\alpha
_{r,d_r}'}$$ En rempla\c cant $\alpha _{i,j}'$ par son expression
par les $\alpha _{i,j}$, on trouve bien que $2\lambda
_{1,d_1-1}-\lambda _{1,d_1}<0$, parce que c'est le coefficient
devant $\alpha _{1,d_1}$.

Consid\'erons maintenant le cas o\`u $(\pi ,V)$ est temp\'er\'ee.
Alors, le raisonnement pr\'ec\'edent pour $\pi $ de carr\'e
int\'egrable reste valable apr\`es avoir remplac\'e $<0$ par $\leq
0$ au d\'etail pr\`es que certains des $\Delta _{\so ,i}$ peuvent
bien \^etre de type $A_{d_i-1}$. Le raisonnement pr\'ec\'edent ne
permet alors plus de conclure que $\lambda $ est combinaison
lin\'eaire \`a coefficients $\leq 0$ de $\Delta _{\so }$. Mais, on
peut toujours \'echanger les facteurs $GL_{k_1}\times\cdots\times
GL_{k_1}$ et $GL_{k_i}\times\cdots\times GL_{k_i}$ de
$\underline{M}$ pour se ramener \`a $i=1$. Par le proposition {\bf
4.4} (iii) et l'hypoth\`ese, $w_{0,1}\lambda $ est \'egalement une
combinaison lin\'eaire \`a coefficients $\leq 0$ de $\Delta (P)$.
Comme
$$\eqalign{w_{0,1}\lambda =&+\lambda_{1,d_1-1}\alpha _{1,1}' +\cdots +\lambda _{1,1}\alpha
_{1,d_1-1}'-\lambda _{1,d_1}\alpha _{1,d_1}'\cr
&+\lambda_{2,1}\alpha _{2,1}'+\cdots \lambda_{r,d_r}\alpha
_{r,d_r}',\cr}$$ on en d\'eduit que $\lambda _{1,d_1}=0$, apr\`es
avoir remplac\'e les $\alpha _{i,j}'$ par les $\alpha _{i,j}$.
Donc, $\lambda $ est bien une combinaison lin\'eaire de $\Delta
(P)$ \`a coefficients $\leq 0$.

Ceci termine la preuve du th\'eor\`eme. \hfill{\fin 2}

\enddocument

\enddocument